\documentclass{elsart3-1}

\usepackage[francais,english]{babel}
\usepackage[latin1]{inputenc}
\usepackage{amssymb}
\usepackage{amsmath}
\usepackage{wasysym}
\usepackage{stmaryrd}
\usepackage{enumerate}

\newtheorem{theorem}{Theorem}

\newtheorem{e-proposition}[theorem]{Proposition}

\newtheorem{e-definition}[theorem]{Definition\rm}

\newtheorem{theoreme}{Th\'eor\`eme}[section]
\newtheorem{lemme}[theoreme]{Lemme}
\newtheorem{proposition}[theoreme]{Proposition}

\newtheorem{remarque}{\it Remarque}

\renewcommand{\theequation}{\arabic{equation}}
\def\og{\leavevmode\raise.3ex\hbox{$\scriptscriptstyle\langle\!\langle$~}}
\def\fg{\leavevmode\raise.3ex\hbox{~$\!\scriptscriptstyle\,\rangle\!\rangle$}}
\newcommand{\R}[1]{\ensuremath{\mathbb{R}^{#1}}}
\newcommand{\vect}[1]{\mathrm{Vect}\left\{#1\right\}}
\newcommand{\refhypo}[1]{{\bf (H\ref{#1})}}
\newcounter{hypotheses}
\newenvironment{Hypotheses}{\begin{enumerate}[{\bf {(H}1{)}}]\setcounter{enumi}{\thehypotheses}}{\setcounter{hypotheses}{\theenumi}\end{enumerate}}
                      
\newcounter{assumptions}
\newenvironment{Assumptions}{\begin{enumerate}[{\bf {(A}1{)}}]\setcounter{enumi}{\theassumptions}}{\setcounter{assumptions}{\theenumi}\end{enumerate}}

\begin{document}
\sloppy

\begin{frontmatter}

	\selectlanguage{francais}
	\title{Un résultat de consistance pour des SVM fonctionnels par interpolation spline}

	\vspace{-2.6cm}
	\selectlanguage{english}
	\title{A consistency result for functional SVM by spline interpolation}

	\author[nath]{Nathalie Villa}
	\ead{villa@univ-tlse2.fr}
	\author[fab]{Fabrice Rossi}
	\ead{Fabrice.Rossi@inria.fr}

	\address[nath]{\'Equipe GRIMM, Université Toulouse Le Mirail, 5 allées Antonio Machado, 31058 Toulouse cedex 9, France.}
	\address[fab]{Projet AxIS, INRIA-Rocquencourt, Domaine de Voluceau, Rocquencourt, BP 105, 78153 Le Chesnay cedex, France.}

	\begin{abstract}
This note proposes a new methodology for function classification with
Support Vector Machine (SVM). Rather than relying on projection on a truncated
Hilbert basis as in our previous work, we use an implicit spline
interpolation that allows us to compute SVM on the derivatives of the studied
functions. To that end, we propose a kernel defined directly on the
discretizations of the observed functions. We show that this method is
universally consistent.

		\vskip 0.5\baselineskip

		\selectlanguage{francais}
		\noindent{\bf R\'esum\'e}
		\vskip 0.5\baselineskip
		\noindent
Nous proposons dans cette note une nouvelle méthode de discrimination de
données fonctionnelles par Support Vector Machine (SVM). Dans nos travaux
antérieurs, nous nous appuyions sur une projection sur une base hilbertienne
tronquée ; nous proposons ici d'utiliser une interpolation spline implicite,
afin de pouvoir construire un SVM sur les dérivées des fonctions
initiales. Pour cela, nous construisons un noyau qui s'applique directement
sur les discrétisations des observations. Nous montrons la consistance
universelle d'une telle approche. 
	\end{abstract}
\end{frontmatter}

\selectlanguage{english}
\section*{Abridged English version}
We emphasize in \cite{rossi_villa_N2006} the interest of using classical SVM
\cite{vapnik_SLT1998} on the derivatives of the original functions for some
kind of data sets (near infra-red spectrometric curves for example). We
propose here a practical and consistent methodology for using SVM for binary
classifications when the regressor is a smooth function. 

Let $(X,Y)$ be a pair of random variables where $X$ takes its values in the
Sobolev space $\mathcal{H}^m([0,1])=\left\{h\in L^2([0,1]):\
  \forall\,j=1,\ldots,m,\ D^jh \text{ exists (in a weak sense) and }
 D^jh \in L^2([0,1])\right\}$ and $Y \in \{-1,1\}$. We are given $n$
observations of this random pair, $(x_1,y_1),\ldots,(x_n,y_n)$; furthermore,
the $x_i$ ($i=1,\ldots,n$) are not completly known as we are only given a
discretization of them: $\mathbf{x}_i=(x_i(t_1),\ldots,x_i(t_d))^T$. 

The main point of this note is to represent the observations of $X$ by a
$L$-spline interpolation for which the derivatives are implicitly calculated
through the discretization. This $L$-spline interpolation minimizes a penalty
defined by a differential operator $L=D^m+\sum_{j=0}^{m-1}a_j D^j$. This
operator allows us to decompose the space $\mathcal{H}^m$ as
$\mathcal{H}^m=\mathcal{H}_0\oplus \mathcal{H}_1$ where
$\mathcal{H}_0=\textrm{Ker}L$ is a $m$-dimensional Hilbert space and
$\mathcal{H}_1$ is a reproducing kernel Hilbert space (RKHS) with kernel $K$. 
$\mathcal{H}_1$ is defined by $m$ boundary conditions (for all $h\in
\mathcal{H}_1$ and all $j=1,\ldots,m$, $B^jh=0$) and the inner product: for all $u$, $v\in \mathcal{H}_1$, $\langle u,v\rangle_1
=\int_{[0,1]} Lu(t)Lv(t)\textrm{d}t$ (see \cite{besse_ramsay_P1986}
or \cite{berlinet_thomasagnan_RKHSPS2004} for further informations about
RKHS). On the space $\mathcal{H}_1$, the $L$-spline representation of a
discretization is given by the following theorem:
\begin{theorem}[\cite{besse_ramsay_P1986}]
\label{interp_splines_anglais}
Let $x \in \mathcal{H}_1$ be a function known at $t_1,\ldots,t_d$. We assume that
the matrix $\mathbf{K}_d=(K(t_i,t_j))_{i,j=1,\ldots,d}$ is positive definite. Then, there
exists a unique interpolation function $h\in \mathcal{H}_1$ at $t_1,\ldots,
t_d$, such that $\| h\|_1\, \leq\, \| u\|_1$ for any interpolation
function $u\in \mathcal{H}_1$. $h$ is given by:  
	$h=\sum_{i=1}^d c_i K(t_i,.)$, where $c=\mathbf{K}^{-1}_d\mathbf{x}$ and $\mathbf{x}=(x(t_1),\ldots,x(t_d))^T$.

	Moreover, if $h_1$ and $h_2$ are the respective interpolation functions of $x_1$ and $x_2\in \mathcal{H}_1$ defined as above then, $\langle h_1,h_2\rangle_1 = \mathbf{x}_1^T\mathbf{K}_d^{-1} \mathbf{x}_2 = \langle \mathbf{x}_1,\mathbf{x}_2\rangle_{(\R{d},\mathbf{K}_d^{-1})}$, where $(\R{d},\mathbf{K}^{-1}_d)$ is $\R{d}$ with the inner product induced by the matrix $\mathbf{K}_d^{-1}$.
\end{theorem}

Let then, for all $i=1,\ldots,n$, $h_i$ be the $L$-spline interpolating the
observation $x_i$ at $t_1,\ldots,t_d$. Provided that
$\mathbf{K}_d=(K(t_i,t_j))_{i,j=1,\ldots,d}$ is positive definite, we can
construct a SVM on $(h_i)_{i=1,\ldots,n}$ through the discretizations
$(\mathbf{x}_i)_{i=1,\ldots,n}$:

\begin{theorem}
\label{th_svm_anglais}
	Let $G_\gamma^d$ be the gaussian kernel with parameter $\gamma$ on
        \R{d} and $G_\gamma^\infty$ the gaussian kernel with parameter
        $\gamma$ on $L^2([0,1])$ ($G_\gamma(u,v)=e^{-\gamma\|
          u-v\|^2_{\R{d}\textrm{ or }L^2}}$). Then, a SVM on the derivatives of $h_1,\ldots,h_n$ (denoted $\phi_h^{n,d}$) defined by
	\[
	\begin{array}{l}
		\max_\alpha \sum_{i=1}^n \alpha_i - \sum_{i,j=1}^n \alpha_i \alpha_j G_\gamma^\infty(Lh_i,Lh_j)\\
		\textrm{with }\qquad	\sum_{i=1}^n \alpha_i y_i =0,\qquad	0\leq\alpha_i \leq C,\ 1\leq i\leq n,
	\end{array}
	\]
	is equivalent to a SVM on the discretizations $\mathbf{x}_1,\ldots,\mathbf{x}_n$ (denoted $\phi_\mathbf{x}^{n,d}$) :
	\[
	\begin{array}{l}
		\max_\alpha \sum_{i=1}^n \alpha_i - \sum_{i,j=1}^n \alpha_i \alpha_j G_\gamma^d\circ\mathbf{K}_d^{-1/2}(\mathbf{x}_i,\mathbf{x}_j)\\
		\textrm{with }\qquad	\sum_{i=1}^n \alpha_i y_i =0,\qquad	0\leq\alpha_i \leq C,\ 1\leq i\leq n.
	\end{array}
\]
\end{theorem}

Finally, we obtain a consistency result for this model:
\begin{theorem}
\label{th_consist_anglais}
Under the assumptions
\begin{Assumptions}
\item $X$ is a bounded random variable taking its values in $\mathcal{H}_1$,
\item $(\tau_d)_d$ is a sequence of discretization points in $[0,1]$ such that, for all $d\geq1$, $\tau_d=\{t_k\}_{k=1,\ldots,d}$, the matrix $\mathbf{K}_d$ is definite positive and $\mathrm{Span}\{K(t,.),\ t\in \cup_{d\geq1} \tau_d\}$ is dense in $\mathcal{H}_1$,
\item $(C_n^d)_n$ is a sequence such that $C_n^d=\mathcal{O}(n^{1-\beta_d})$ for a $0<\beta_d<1/d$,
\end{Assumptions}
The sequence of SVM classifiers $\phi_h^{n,d}$ defined as in Theorem~\ref{interp_splines_anglais}, with $C=(C_n^d)_n$, is universally consistant in \R{d}, that is:
\[
\lim_{d\rightarrow +\infty} \lim_{n\rightarrow +\infty}Err\phi_h^{n,d}=Err^*
\]
where $Err^*$ is the Bayes error, $\inf_{\phi:\mathcal{H}_1\rightarrow\{-1,1\}}\mathbb{P}(\phi(X)\neq Y)$, and $Err\phi$ is the error of a classifier $\phi$, $\mathbb{P}(\phi(X)\neq Y)$.
\end{theorem}

\selectlanguage{francais}
\section{Introduction}

Nous nous intéressons ici à l'utilisation des SVM pour le traitement de
données fonctionnelles. De manière plus précise, il s'agit de résoudre des
problèmes de discrimination binaire pour lesquels la variable explicative est
fonctionnelle. Nous montrons dans \cite{rossi_villa_N2006} l'intérêt pratique,
pour certains types de données, d'utiliser des SVM (Support Vector Machine,
voir \cite{vapnik_SLT1998}) sur les dérivées des fonctions initiales ; nous
proposons, dans cette note, une méthodologie permettant de mettre en \oe uvre
un tel traitement et démontrons un résultat de consistance universel associé à
celle-ci. 

Pour cela, nous étudions un couple de variables aléatoires $(X,Y)$ où $X$ est
supposée \og régulière\fg\ et prend ses valeurs dans l'espace de Sobolev
$\mathcal{H}^m([0,1])=\left\{h\in L^2([0,1])\ :\ \forall\,j=1,\ldots,m,\ D^jh
\right.$ existe (au sens faible) et $\left. D^jh \in L^2([0,1])\right\}$ et $Y
\in \{-1,1\}$. Ce couple est connu grâce à $n$ observations,
$(x_1,y_1),\ldots,(x_n,y_n)$ ; en fait, les $x_i$ ($i=1,\ldots,n$) ne sont pas
connues de manière exacte mais simplement au travers d'une discrétisation
$\mathbf{x}_i=(x_i(t_1),\ldots,x_i(t_d))^T$ (les points de discrétisation sont
les mêmes pour tous les $x_i$ et sont déterministes). Le problème est alors de
construire, à partir de ces données, un classifieur capable de prédire $Y$
connaissant $X$.  En tirant partie de la structure d'espace de Hilbert à noyau
reproduisant (RKHS) de $\mathcal{H}^m([0,1])$, les observations de $X$ seront
représentées par une interpolation spline sur laquelle les dérivées
s'expriment de manière naturelle en fonction de la discrétisation.

\section{Interpolation $L$-Spline}

On choisit de représenter les observations de $\mathcal{H}^m([0,1])$ à travers une interpolation $L$-spline : celle-ci interpole exactement la fonction aux points de discrétisation tout en minimisant une pénalité définie à partir d'un opérateur différentiel $L=D^m +\sum_{j=0}^{m-1} a_j D^j$. On peut montrer que, si le noyau de cet opérateur, $\textrm{Ker}L=\mathcal{H}_0$ est un sous-espace de dimension $m$ de $\mathcal{H}^m$, on peut écrire $\mathcal{H}^m=\mathcal{H}_0 \oplus \mathcal{H}_1$ où $\mathcal{H}_1$ est un sous-espace vectoriel de $\mathcal{H}^m$ défini par $m$ conditions aux bornes, $\forall\,h\in\mathcal{H}_1$ et $\forall\,j=1,\ldots,m$, $B^jh=0$, et muni du produit scalaire $\forall\,u,v\in \mathcal{H}_1$, $\langle u,v\rangle_1=\langle Lu,Lv\rangle_{L^2}=\int_{[0,1]} Lu(t) Lv(t)\,\textrm{d}t$ (voir, par exemple, \cite{besse_ramsay_P1986} ou \cite{berlinet_thomasagnan_RKHSPS2004}). $\mathcal{H}_0$ et $\mathcal{H}_1$ sont deux espaces de Hilbert à noyau reproduisant et on note $K$ le noyau reproduisant de $\mathcal{H}_1$ ; on donne, dans \cite{villa_rossi_JDS2006}, des exemples de décompositions de $\mathcal{H}^m$ et on explique, sur ces exemples, comment calculer $K$.

Cette décomposition permet de définir simplement le produit scalaire entre les
représentations des fonctions à partir des discrétisations initiales :
\begin{theoreme}[\cite{besse_ramsay_P1986}]
\label{interp_splines}
	Soit $x\in \mathcal{H}_1$ une fonction connue aux points de
        discrétisation $t_1,\ldots, t_d$. Supposons, en outre, que la matrice
        $\mathbf{K}_d=(K(t_i,t_j))_{i,j}$ soit définie positive. Alors, il
        existe une unique fonction d'interpolation $h\in \mathcal{H}_1$ aux
        points $t_1,\ldots, t_d$ telle que $\| h\|_1\, \leq\, \| u\|_1$ pour
        toute fonction d'interpolation $u\in \mathcal{H}_1$. $h$ est donnée
        par : 
	\[
	h=\sum_{i=1}^d c_i K(t_i,.)
	\]
	où $c=\mathbf{K}^{-1}_d\mathbf{x}$ avec $\mathbf{x}=(x(t_1),\ldots,x(t_d))^T$.

	De plus, si $h_1$ et $h_2$ sont les deux fonctions d'interpolation de $x_1$ et $x_2\in \mathcal{H}_1$ comme définies ci-dessus, alors
	\begin{equation}
	\label{ps_dans_H}
		\langle h_1,h_2\rangle_1 = \mathbf{x}_1^T\mathbf{K}^{-1}_d \mathbf{x}_2 = \langle \mathbf{x}_1,\mathbf{x}_2\rangle_{(\R{d},\mathbf{K}^{-1}_d)}
	\end{equation}
	où $(\R{d},\mathbf{K}^{-1}_d)$ est l'espace $\R{d}$ muni du produit scalaire induit par la matrice $\mathbf{K}^{-1}_d$.
\end{theoreme}

La fonction d'interpolation spline est donc simplement
$h=\mathcal{P}_{\vect{K(t_k,.),\ k=1,\ldots,d}}(x)$, où $\mathcal{P}_V$ est
l'opérateur de projection orthogonale sur $V$ dans $\mathcal{H}_1$, ce qui rapproche la
méthodologie proposée ici de celle développée dans \cite{rossi_villa_N2006} et
inspirée des travaux de \cite{biau_bunea_webkamp_IEEETIT2005}. 
Ceci permet de déterminer la perte d'information induite par l'interpolation,
notamment en terme de perturbation de l'erreur de Bayes, comme le montre le
résultat suivant :

\begin{lemme}
\label{lemme1}
	Soient 
	\begin{Hypotheses}
		\item $X$ une variable aléatoire à valeurs dans $\mathcal{H}_1$ ;
		\item \label{suite_reg} $(\tau_d)_{d\geq1}$ une suite de points de discrétisation de $[0,1]$ telle que $\forall\,d\geq1$, $\tau_d=\{t_k\}_{k=1,\ldots,d}$, la matrice $\mathbf{K}_d=(K(t_i,t_j))_{i,j=1,\ldots,d}$ est inversible et $\vect{K(t,.),\ t\in \cup_{d\geq 1}\tau_d}$ est dense dans $\mathcal{H}_1$.
	\end{Hypotheses}
On note $V_d=\vect{K(t,.),\ t\in\tau_d}$ et
$\mathcal{P}_d(x)=\mathcal{P}_{V_d}(x)$.  On a alors
	\begin{equation}
	\label{ineg1}
		\lim_{d\rightarrow +\infty} Err_d^* = Err^*
	\end{equation}
avec $Err_d^*=\inf_{\phi:V_d\rightarrow \{-1,1\}}
\mathbb{P}(\phi(\mathcal{P}_d(X))\neq Y)$ (erreur de Bayes de la
représentation $L$-spline), et $Err^*$ est l'erreur de Bayes donnée par :
$\inf_{\phi:\mathcal{H}_1\rightarrow\{-1,1\}}\mathbb{P}(\phi(X)\neq Y)$. 
\end{lemme}

\noindent {\it Démonstration :}  Les $\vect{K(t,.),\ t\in\tau_d}$ ($d\geq1$) sont des ensembles emboîtés et, par densité, $\forall\, x\in \mathcal{H}_1$, $\lim_{d\rightarrow+\infty} \mathcal{P}_d(x) = x$ dans $\mathcal{H}_1$.

Par ailleurs, les $\sigma$-algèbres
$\sigma(\mathcal{P}_d(X))=\sigma(\mathbf{K}^{-1}_d (X(t_1),\ldots, X(t_d))^T)$
forment clairement une filtration. Comme $\mathbb{E}(|Y|)\leq 1$,
$\mathbb{E}(Y|\mathcal{P}_d(X))$ est une martingale uniformément intégrable
pour cette filtration (cf \cite{pollard_UGMTP2002} lemme 35 page 154), cette
martingale converge en norme $L^1$ vers $\mathbb{E}(Y|\sigma(\cup_d
\sigma (\mathcal{P}_d(X))))$ (cf théorème 36 page 154 de \cite{pollard_UGMTP2002}),
dont la valeur est $\mathbb{E}(Y|X)$ (puisque $\mathcal{P}_d(X)$ est fonction de $X$, $\sigma(\cup_d
\sigma (\mathcal{P}_d(X)))\subset \sigma(X)$ et, inversement, $X$ est $\sigma(\cup_d
\sigma (\mathcal{P}_d(X)))$-mesurable comme limite des variables aléatoires $(\mathcal{P}_d(X))_d$, $\sigma(\cup_d
\sigma (\mathcal{P}_d(X)))$-mesurables).

Nous concluons en utilisant l'inégalité classique $Err_d^* - Err^* \leq
2\mathbb{E}|\mathbb{E}(Y|\mathcal{P}_d(X))-\mathbb{E}(Y|X)|$ (cf e.g.
\cite{devroye_gyorfi_lugosi_PTPR1996}, théorème 2.2).
$\Square$

\section{SVM sur dérivées}

Notons, $\forall\,i=1,\ldots,n$, $h_i$ la spline d'interpolation de l'observation $x_i$ aux points de discrétisation $t_1,\ldots,t_d$ définie comme dans le Théorème~\ref{interp_splines}. Alors, si la matrice $\mathbf{K}_d=(K(t_i,t_j))_{i,j=1,\ldots,d}$ est inversible, on peut définir un SVM sur les dérivées des $L$-splines d'interpolation par le théorème suivant :
\begin{theoreme}
\label{th_svm}
	Soit $G_\gamma^d$ le noyau gaussien de paramètre $\gamma$ sur \R{d} et $G_\gamma^\infty$ le noyau gaussien de paramètre $\gamma$ sur $L^2([0,1])$ ($G_\gamma(u,v)=e^{-\gamma\| u-v\|^2_{\R{d}\textrm{ ou }L^2}}$). Alors, le SVM sur les dérivées des fonctions $h_1,\ldots,h_n$ (noté $\phi_h^{n,d}$) défini par
	\[
	\begin{array}{l}
		\max_\alpha \sum_{i=1}^n \alpha_i - \sum_{i,j=1}^n \alpha_i \alpha_j G_\gamma^\infty(Lh_i,Lh_j)\\
		\textrm{avec }\qquad	\sum_{i=1}^n \alpha_i y_i =0,\qquad	0\leq\alpha_i \leq C,\ 1\leq i\leq n,
	\end{array}
	\]
est équivalent au SVM sur les discrétisations $\mathbf{x}_1,\ldots,\mathbf{x}_n$ (noté $\phi_\mathbf{x}^{n,d}$) :
	\[
	\begin{array}{l}
		\max_\alpha \sum_{i=1}^n \alpha_i - \sum_{i,j=1}^n \alpha_i \alpha_j G_\gamma^d\circ\mathbf{K}_d^{-1/2}(\mathbf{x}_i,\mathbf{x}_j)\\
		\textrm{avec }\qquad	\sum_{i=1}^n \alpha_i y_i =0,\qquad	0\leq\alpha_i \leq C,\ 1\leq i\leq n.
	\end{array}
\]
\end{theoreme}

\noindent{\it Démonstration :} Il suffit de constater, d'après (\ref{ps_dans_H}), que $\forall\,i,j=1,\ldots,n$, $G_\gamma^\infty(Lh_i,Lh_j)=e^{-\gamma\| Lh_i-Lh_j\|^2_{L^2}}=e^{-\gamma\| \mathbf{x}_i-\mathbf{x}_j\|^2_{(\R{d},\mathbf{K}^{-1})}}=e^{-\gamma\| \mathbf{K}_d^{-1/2}\mathbf{x}_i-\mathbf{K}_d^{-1/2}\mathbf{x}_j\|^2_{\R{d}}}$.$\Square$

Or, \cite{steinwart_JC2002} démontre la consistance universelle des SVM $d$-dimensionnels. Ainsi, à suite de discrétisation fixée $t_1,\ldots,t_d$, on peut démontrer la consistance universelle des SVM $\phi_h^{n,d}$ vers l'erreur de Bayes de la représentation $L$-spline ; ainsi, à discrétisation fixée, $\phi_h^{n,d}$ est asymptotiquement optimal :
\begin{lemme}
\label{lemme2}
	Soit $t_1,\ldots,t_d$ des points de discrétisation tels que $\mathbf{K}_d=(K(t_i,t_j))_{i,j=1,\ldots,d}$ est inversible. Supposons que
	\begin{Hypotheses}
		\item $(C_n^d)_n$ est une suite telle que $C_n^d=\mathcal{O}(n^{1-\beta_d})$ pour $0<\beta_d<1/d$ ;
		\item $X$ est une variable aléatoire bornée dans $\mathcal{H}_1$.
	\end{Hypotheses}
	Alors, le SVM $\phi_h^{n,d}$ défini comme dans le Théorème~\ref{th_svm},
avec pour paramètre $C=C_n^d$, est universellement consistant
dans \R{d} :
	\begin{equation}
	\label{lim_etape2}
		\lim_{n\rightarrow+\infty} Err\phi_h^{n,d}=Err^*_d
	\end{equation}
	pour $Err\phi=\mathbb{P}(\phi(X)\neq Y)$.
\end{lemme}

\noindent {\it Démonstration :} On note
$\mathbf{X}=(X(t_1),\ldots,X(t_d))^T$. Par le Théorème \ref{th_svm}, $Err\phi_h^{n,d}=Err\phi_\mathbf{x}^{n,d}$ et, puisque $\mathbf{K}_d$ est inversible, $\inf_{\phi:\R{d}\rightarrow \{-1,1\}} \mathbb{P}(\phi(\mathbf{X})\neq Y)=\inf_{\phi:V_d\rightarrow \{-1,1\}} \mathbb{P}(\phi(\mathcal{P}_d(X))\neq Y)=Err^*_d$. Or, d'après \cite{steinwart_JC2002}, les SVM dans \R{d} sont universellement consistants ; pour cela, on doit vérifier :
\begin{enumerate}[1{.}]
	\item la variable aléatoire explicative prend ses valeurs dans un
          compact de \R{d} : comme $X$ prend ses valeurs dans un borné de
          $\mathcal{H}_1$, $\mathbf{X}$ prend ses valeurs dans un borné de
          $\R{d}$, c'est-à-dire, un compact de \R{d}, noté $U$ ;
\item le noyau utilisé doit être \emph{universel} : Steinwart montre dans
  \cite{steinwart_JMLR2001} que le noyau gaussien $d$-dimensionnel est
  universel. Il montre aussi que tout noyau obtenu en composant une fonction
  continue et injective avec un noyau universel est lui aussi universel. Or,
  $\mathbf{K}_d^{-1/2}$ est continue et injective, et donc le noyau
  $G_\gamma^d\circ \mathbf{K}_d^{-1/2}$ est universel : l'ensemble des
  fonctions de la forme $\langle \Phi\circ
  \mathbf{K}^{-1/2}_d(.),w\rangle_\mathcal{X}$ ($w\in \mathcal{X}$) est dense
  dans l'ensemble des fonctions continues sur un compact de \R{d} (où
  $\mathcal{X}$ désigne le RKHS associé au noyau $G_\gamma^d\circ
  \mathbf{K}_d^{-1/2}$) ; 
	\item on doit contrôler le nombre de couverture
          $\mathcal{N}(G_\gamma^d\circ \mathbf{K}^{-1/2}_d,\epsilon)$,
          c'est-à-dire le nombre minimal de boules de rayon $\epsilon$ (au sens
          de la métrique de \R{d} définie par le noyau $G_\gamma^d\circ
          \mathbf{K}^{-1/2}_d$) nécessaires pour recouvrir $U$ le support compact de $\mathbf{X}$. Or, on montre aisément que $\mathcal{N}(G_\gamma^d\circ
          \mathbf{K}^{-1/2}_d,\epsilon)\leq
          \mathcal{N}(G_\gamma^d,\epsilon)$, puis on utilise
          \cite{steinwart_JC2002} pour obtenir $\mathcal{N}(G_\gamma^d,
          \epsilon)=\mathcal{O}_n(\epsilon^{-d})$ et donc 
          $\mathcal{N}(G_\gamma^d\circ
          \mathbf{K}^{-1/2}_d,\epsilon)=\mathcal{O}_n(\epsilon^{-d})$ ;
	\item la suite $(C_n^d)_n$ est bien de la forme requise ($\mathcal{O}(n^{1-\beta_d})$ avec $0<\beta_d<1/d$).
\end{enumerate}
On conclut donc, par le Théorème~2 de \cite{steinwart_JC2002}, que $Err\phi_h^{n,d}=Err\phi_\mathbf{x}^{n,d}\xrightarrow{n\rightarrow+\infty} \inf_{\phi:\R{d}\rightarrow \{-1,1\}} \mathbb{P}(\phi(\mathbf{X})\neq Y) = Err_d^*$.$\Square$

\section{Consistance}
L'utilisation de noyaux définis comme dans le Théorème~\ref{th_svm} sous les hypothèses formulées dans les lemmes \ref{lemme1} et \ref{lemme2} conduit à des SVM universellement consistants (double limite lorsque le nombre de points de discrétisation tend vers l'infini et le nombre d'observations tend vers l'infini) :
\begin{theoreme}
\label{th_consist}
	Sous les hypothèses {\bf (H1)-(H4)}, le SVM défini comme dans le Théorème~\ref{th_svm}, $\phi_h^{n,d}$, pour les points d'interpolation $(\tau_d)_{d\geq 1}$ et la suite $C=(C_n^d)_n$ est universellement consistant dans $\mathcal{H}_1$ :
	\[
	\lim_{d\rightarrow+\infty} \lim_{n\rightarrow +\infty} Err \phi_h^{n,d} =Err^*.
	\]
\end{theoreme}

\noindent{\it Démonstration :} On écrit $Err\phi_h^{n,d}-Err^*=(Err\phi_h^{n,d}-Err_d^*)+(Err_d^*-Err^*)$. Soit alors $\epsilon >0$. Par le Lemme~\ref{lemme1}, il existe $D_0>0$ : $\forall\,d\geq D_0$, $Err_d^*-Err^*\leq \epsilon$. Soit alors $d\geq D_0$ ; par le Lemme~\ref{lemme2}, $\exists N_0>0$ : $\forall\,n\geq N_0$, $(Err\phi_h^{n,d}-Err_d^*)\leq \epsilon$, ce qui conclut la preuve.$\Square$

\begin{remarque}
	La discrétisation des fonctions est en général induite par le problème. Si $\tau$ est une discrétisation donnée, on peut supposer, quitte à retirer quelques points, que la matrice $(K(t,t'))_{t,t'\in \tau}$ est inversible. Il existe alors une suite de points de discrétisation telle que $\tau=\tau_1$ et qui vérifie l'hypothèse \refhypo{suite_reg} :
\end{remarque}

\begin{proposition}
	Si $\tau$ est un ensemble fini de points de $[0,1]$ tels que $(K(t,t'))_{t,t'\in \tau}$ est inversible alors, il existe un ensemble dénombrable $\mathcal{D}_0=(t_k)_{k\geq 1}\subset[0,1]$ tel que
	\begin{itemize}
		\item $\tau \subset \mathcal{D}_0$ ;
		\item $\vect{K(t,.),\ t\in\mathcal{D}_0}$ est dense dans $\mathcal{H}_1$ ;
		\item pour tout $d\geq 1$, la matrice $(K(t_i,t_j))_{i,j=1,\ldots,d}$ est inversible.
	\end{itemize}
\end{proposition}

\noindent {\it Démonstration :} Par le Théorème~15 de \cite{berlinet_thomasagnan_RKHSPS2004}, l'espace de Hilbert $\mathcal{H}_1$ est séparable (comme ensemble de fonctions continues) dès que $m\geq1$. Or, $(K(t,t'))_{t,t'\in \tau}$ est inversible est équivalent au fait que $\{K(t,.),\ t\in \tau\}$ est une famille de fonctions linéairement indépendantes. Ainsi, par le Théorème~8 de \cite{berlinet_thomasagnan_RKHSPS2004}, il existe un support dénombrable de $\mathcal{H}$ contenant $\tau$, c'est-à-dire, un ensemble dénombrable $\mathcal{D}_0$ tel que $\tau\subset \mathcal{D}_0$, les $\{K(t,.),\ t\in\mathcal{D}_0\}$ sont linéairement indépendants et $\vect{K(t,.),\ t\in \mathcal{D}_0}$ est dense dans $\mathcal{H}_1$.$\Square$

\begin{remarque}
	En pratique, la matrice $(K(t,t'))_{t,t'\in \tau}$ est souvent mal conditionnée dès que le cardinal de $\tau$ est élevé. Ainsi, il sera donc préférable d'introduire un paramètre de régularisation (splines de lissage) afin de permettre l'inversion de celle-ci.
\end{remarque}


\section*{Remerciements}
Les auteurs tiennent à remercier les deux rapporteurs pour leurs recommandations pertinentes qui ont permis l'amélioration de cette note.



\end{document}